\documentclass[12pt]{amsart}

\usepackage{amssymb}

\usepackage{enumitem}

\usepackage{graphicx}

\makeatletter
\@namedef{subjclassname@2020}{%
  \textup{2020} Mathematics Subject Classification}
\makeatother

\usepackage[T1]{fontenc}

\newtheorem{theorem}{Theorem}[section]
\newtheorem{corollary}[theorem]{Corollary}

\newtheorem{proposition}[theorem]{Proposition}

\theoremstyle{definition}
\newtheorem{definition}[theorem]{Definition}
\newtheorem{remark}[theorem]{Remark}
\newtheorem{example}[theorem]{Example}

\numberwithin{equation}{section}

\begin{document}


\baselineskip=17pt


\title[Strong Pseudoconvexity in Banach Spaces]{Strong Pseudoconvexity in Banach Spaces}

\author{Sofia Ortega Castillo}

\thanks{This author was under a CONACYT fellowship during part of the writing period of this manuscript. Additionally, the author was supported by a UUKi Rutherford Strategic Partner Grant to do a research visit at the University of Bath, and would also like to thank the Erwin Schr\"{o}dinger International Institute for Mathematics and Physics for granting them support during a workshop visit.}

\address{Centro de Investigaci\'on en Matem\'aticas A.~C., Jalisco S/N, Colonia Valenciana, Guanajuato, 36023, Guanajuato, M\'exico}

\address{Department of Mathematical Sciences, University of Bath, Claverton Down, Bath, BA2 7AY, United Kingdom}
            
\address{Departamento de Matem\'aticas, Centro Universitario de Ciencias Exactas e Ingenier\'ias, Universidad de Guadalajara, Blvd. Marcelino Garc\'ia Barrag\'an \#1421, Esq. Calzada Ol\'impica, Guadalajara, 44430, Jalisco, M\'exico}

\email{sofia.ortega@academicos.udg.mx}



\date{}

\begin{abstract}
Having been unclear how to define that a domain is strictly pseudoconvex in the infinite-dimensional setting, we develop a general theory having Banach spaces in mind. We first focus on finite dimension and eliminate the need of two degrees of differentiability of the boundary of a domain, since differentiable functions are difficult to find in infinite dimension. We introduce $\ell$-strict pseudoconvexity for $\ell\geq 1$,  $1$-strict pseudoconvexity at the boundary, $\ell$-uniform pseudoconvexity for $\ell\geq 0$ and finally strong pseudoconvexity. Defining $\ell$-strict pseudoconvexity and $\ell$-uniform pseudoconvexity for $\ell<2$ depends on extending a notion of strict plurisubharmonicity to cases lacking $C^2$-smoothness, first studying it in the sense of distribution and then considering it in infinite dimension. Examples of strictly plurisubharmonic functions as well as strongly pseudoconvex domains are presented, which end up related to important classical Banach spaces. Finally, some solutions to the inhomogeneous Cauchy-Riemann equations for $\overline{\partial}$-closed $(0,1)$-forms in infinite-dimensional domains are shown, giving new information about domains affinely isomorphic to the ball of $\ell_1$ which appeared in the study of strong pseudoconvexity and some more domains biholomorphically related to open and convex domains of $\ell_1$ such as its ball.
\end{abstract}

\subjclass[2020]{Primary 32T15, 46B25; Secondary 46F10, 31C10}

\keywords{strictly plurisubharmonic on average continuously, $\ell$-strictly pseudoconvex, $\ell$-uniformly pseudoconvex, strongly pseudoconvex}

\maketitle

\section{Introduction}

Classically, strongly pseudoconvex domains provide concrete examples of domains of holomorphy in several complex variables that have additional structure. Domains of holomorphy are  those domains $\Omega$ whose boundary points are each a singularity for a holomorphic function on $\Omega$ \cite{M}. Basic examples of domains of holomorphy include any domain in the complex plane, domains of convergence of multivariable power series and any convex domain in a Banach space. Meanwhile, the well known Hartogs' domain is not a domain of holomorphy because it allows analytic continuation to a strictly larger domain \cite[Ch.~II, \S 1.1]{R}.

\medskip

In the Euclidean complex space $\mathbb{C}^n$, when the boundary of a given domain $U$ has two degrees of smoothness we have that $U$ being a domain of holomorphy is characterized by a complex-differential property at boundary points which corresponds to a complex analogue of a well-known differential condition satisfied by convex domains, which is usually introduced as (Levi) pseudoconvexity. Thus strong pseudoconvexity is presented using the complex analogue of the differential condition that defines strict convexity \cite[Ch.~II, \S 2.6]{R}. A simpler equivalent condition, that reduces to the strict plurisubharmonicity of a $C^2$ function defined in a neighborhood of the boundary, has become common too \cite[Ch.~II, \S 2.8]{R}, \cite[\S 1.1]{Ke}. Strictly plurisubharmonic functions are those $C^2$ functions whose complex Hessian is positive definite, where the complex Hessian is a complex analogue of the Hessian which is also called the Levi form.

\medskip

Still in finite dimension, pseudoconvexity was later extended to domains without $C^2$ boundary as admitting a plurisubharmonic exhaustion function on the whole domain \cite[Ch.~II, \S 5.4]{R}, where an exhaustion function is one that has relatively compact sublevel sets. Once in an infinite-dimensional Banach space $X$, the notion of pseudoconvexity of an open set $U$ has been extended by the plurisubharmonicity of $-log d_U$, where $d_U$ denotes the distance to the boundary of $U$. A list of other equivalent ways to define pseudoconvexity in infinite dimension is in \cite[Ch.~VIII, \S 37]{M}. In particular, pseudoconvexity can be characterized by its behavior on finite-dimensional spaces, that is, $U$ is pseudoconvex if and only if $U\cap M$ is pseudoconvex for each finite-dimensional subspace $M$ of $X$. In infinite dimension still every domain of holomorphy is pseudoconvex, but there is a nonseparable Banach space for which the converse is false, while for general separable Banach spaces this problem remains open \cite[Ch.~VIII, \S 37]{M}. In separable Banach spaces with the bounded approximation property, indeed pseudoconvex domains are domains of holomorphy \cite[Ch.~X, \S 45]{M}.

\medskip

In this article we look at a generalization of strong pseudoconvexity to the infinite dimensional setting, for which we first consider related notions in the case that the boundary of a given domain does not have two degrees of smoothness. We were initially interested in determining whether the unit ball of $\ell_1$ has enough structure to be considered strongly pseudoconvex, since it was shown in \cite{L} that such ball admits solutions to the (inhomogeneous) Cauchy-Riemann equations with Lipschitz conditions although the Cauchy-Riemann equations with simply bounded conditions are not solvable, not even locally. Meanwhile, the Cauchy-Riemann equations with bounded conditions on certain finite-dimensional strongly pseudoconvex domains admit solutions that extend continuously to the boundary \cite{Ke}. Having solutions that extend continuously to the boundary to the bounded Cauchy-Riemann equations is a helpful tool to study the boundary behavior of bounded holomorphic functions on a ball of a Banach space \cite{Mc}. And the study of the boundary behavior of bounded holomorphic functions on the ball of a separable Banach spaces in a sense reduces to studying the boundary behavior of bounded holomorphic functions on the ball of Banach spaces that are $\ell_1$-sums of finite-dimensional spaces \cite{JO}. Thus we focus on proving the strong pseudoconvexity of $B_{\ell_1}$ and $B_{\ell_1^n}$, yet we will see that further examples can be obtained as the image of affine isomorphisms of $B_{\ell_1}$ and $B_{\ell_1^n}$, and that we can solve the Cauchy-Riemann equations for $\overline{\partial}$-closed $(0,1)$-forms on such domains. Other examples of strongly pseudoconvex domains exhibited in this article include $B_{\ell_p}$ and $B_{\ell_p^n}$ ($1<p \leq 2$, $n\in \mathbb{N}$). It remains an open problem to find further examples of domains strongly pseudoconvex and to find out whether we can solve the Cauchy-Riemann equations for $\overline{\partial}$-closed $(0,1)$-forms on arbitrary strongly pseudoconvex domains in a Banach space, such as $B_{\ell_p}$ and $B_{\ell_p^n}$ ($1<p \leq 2$, $n\in \mathbb{N}$).

\medskip

Foundations on plurisubharmonicity and pseudoconvexity in Banach spaces, as well as a basic treatment of distributions, can be found in \cite{M}. The reader interested in a deep study of pseudoconvexity in $\mathbb{C}^n$ will find it in \cite{Si}.

\section{Strict plurisubharmonicity}

From now on, let $X$ denote a complex Banach space with open unit ball $B_X$ and norm $\|\cdot\|$, let $U$ denote an open subset of $X$ with boundary $bU$, and let $d_U$ denote the distance function to $bU$. We will also denote with $m$ the Lebesgue measure in $\mathbb{C}^n$ seen as $\mathbb{R}^{2n}$.

\begin{definition}
A function $f:U\to[-\infty,\infty)$ is called \textit{plurisubharmonic} if $f$ is upper semicontinuous and for each $a\in U$ and $b\in X$ such that $a+\overline{\mathbb{D}}\cdot b\subset U$ we have that
$$f(a)\leq \frac{1}{2\pi} \int_0^{2\pi} f(a+e^{i\theta}b) d\theta.$$
\end{definition}

Given a differentiable mapping $f:U\to\mathbb{R}$ and $a\in U$, we will write $Df(a)$ for the Fr\'echet derivative of $f$ at $a$, and in turn its complex-linear and complex-antilinear parts will be denoted by $D'f(a)$ and $D''f(a)$, respectively, which are given by
$$D'f(a)(b)=1/2[Df(a)(b)-iDf(a)(ib)],$$
$$D''f(a)(b)=1/2[Df(a)(b)+iDf(a)(ib)],$$
for every $b \in X$.

\medskip

It is known that a function $f \in C^2(U,\mathbb{R})$ is plurisubharmonic if and only if its complex Hessian is positive semi-definite, i.e. for each $a \in U$ and $b\in X$ we have that
\begin{equation}\label{eq2.1}
D'D''f(a)(b,b)\geq 0.
\end{equation}

\begin{definition}
A function $f\in C^2(U,\mathbb{R})$ is called \textit{strictly plurisubharmonic} when the complex Hessian of $f$ is positive definite, i.e. when a proper inequality in (\ref{eq2.1}) for $b\neq 0$ is satisfied \cite[\S 35]{M}.
\end{definition}

When we are in the Euclidean complex space $\mathbb{C}^n$, we aim to understand a suitable extension of strict plurisubharmonicity to distributions, using that in finite dimension a function $f\in C^2(U)$ is strictly plurisubharmonic if and only if there exists $ \psi \in C(U)$ positive such that
$$D'D''f(a)(b,b)\geq \psi(a) \|b\|^2 \mbox{ for all } a\in U \mbox{ and }b\in \mathbb{C}^n.$$

\medskip

Due to this observation, we introduce the following notion.

\begin{definition}
In the arbitrary Banach space setting we will say that $f\in C^2(U,\mathbb{R})$ is \textit{strictly plurisubharmonic continuously} when there exists  $\psi \in C(U)$ positive such that $D'D''f(a)(b,b)\geq \psi(a) \|b\|^2 \mbox{ for all } a\in U \mbox{ and }b\in X.$
\end{definition}

If $U$ is an open subset of $\mathbb{C}^n$, we will denote the real-valued test functions on $U$ by $\mathcal{D}(U)$. A distribution on $U$ is known to be a continuous functional on $\mathcal{D}(U)$. We shall denote by $\mathcal{D}'(U)$ the vector space of all distributions on $U$.

\begin{definition}
If $U\subset \mathbb{C}^n$, given $f \in L^1(U,\textit{loc})$, we say that $f$ is (strictly) plurisubharmonic in distribution if the distribution it induces is (strictly) plurisubharmonic. At the same time, a distribution $T \in \mathcal{D}'(U)$ is called plurisubharmonic if
$$D'D''T(\phi)(w,w):=\sum_{j,k=1}^{n} \frac{\partial^2 T}{\partial z_j \partial \overline{z_k}}(\phi)w_j \overline{w_k}\geq 0, \mbox{ for all } \phi\geq 0 \mbox{ in } \mathcal{D}(U) \mbox{ and } w \in \mathbb{C}^n.$$

And we will say that $T \in \mathcal{D}'(U)$ is strictly plurisubharmonic if there exists $ \psi \in C(U)$ positive such that 
\begin{equation}\label{eq2}
D'D''T(\phi)(w,w)\geq \Big(\int_U \psi \cdot \phi \; dm\Big) \|w\|^2, \mbox{ for all } \phi\geq 0 \mbox{ in } \mathcal{D}(U) \mbox{ and } w \in \mathbb{C}^n.
\end{equation}
\end{definition}

It has been proved, e.~g. in \cite[\S 3.2 and 4.1]{H2}, that plurisubharmonicity is equivalent to plurisubharmonicity in distribution in the following sense: 

\medskip

Suppose that U is a connected domain in $\mathbb{C}^n$. If $f \neq -\infty$ is plurisubharmonic on $U$, then $f \in L^1(U,\textit{loc})$ and $f$ is plurisubharmonic in distribution. Conversely, if $T \in \mathcal{D}'(U)$ is plurisubharmonic then there exists $f \in L^1(U,\textit{loc})$ plurisubharmonic such that $f$ induces the distribution $T$. As a corollary, if $f \in L^1(U,\textit{loc})$ is plurisubharmonic in distribution then there exists $g \in L^1(U,\textit{loc})$ plurisubharmonic such that $f=g$ $m$-a.e.

\medskip

To prove an analogous version of such result for strict plurisubharmonicity, let us introduce the following concept.

\begin{definition}
We will say that an upper semicontinuous function $g: U\subset X \to [-\infty,\infty)$ is \textit{strictly plurisubharmonic on average (continuously)} if there exists $\varphi \in C(U)$ positive such that for all $a \in U$ and $b \in X$ of small norm (with size depending lower semicontinuously on $a$) we have that,
\begin{equation}\label{eq4.1}
\varphi(a)\|b\|^2 +g(a) \leq  \frac{1}{2\pi}\int_0^{2\pi} g(a+e^{i\theta}b)d\theta.
\end{equation}
If $g$ is as above, let us specifically call it \textit{strictly plurisubharmonic on average (continuously) with respect to the function $\varphi$}.
\end{definition}

\begin{proposition}\label{thm4.1}
Suppose that U is a bounded connected domain in $\mathbb{C}^n$. If $f:U\to [-\infty,\infty)$, with $f\neq -\infty$, is strictly plurisubharmonic on average then $f \in L^1(U,\textit{loc})$ and $f$ is strictly plurisubharmonic in distribution. Conversely, if $T \in \mathcal{D}'(U)$ is strictly plurisubharmonic then there exists $f \in L^1(U,\textit{loc})$ strictly plurisubharmonic on average such that $f$ induces the distribution $T$. As a consequence, if $f \in L^1(U,\textit{loc})$ is strictly plurisubharmonic in distribution then there exists $g \in L^1(U,\textit{loc})$ strictly plurisubharmonic on average such that $f=g$ $m$-a.e.
\end{proposition}
\begin{proof}

If $f\neq -\infty$ is strictly plurisubharmonic on average, then $f$ is in particular plurisubharmonic, so we can use the relationship to plurisubharmonicity in distribution to deduce that $f \in L^1(U,\textit{loc})$. Moreover, since $f$ is strictly plurisubharmonic on average in $U\subset \mathbb{C}^n$, there exists a positive function $\psi \in C(U)$ such that
$$\psi(a)\|b\|^2+f(a) \leq \frac{1}{2\pi}\int_0^{2\pi} f(a+e^{i\theta}b)d\theta,$$
for all $a \in U \mbox{ and } b \in \mathbb{C}^n$ of norm less than $\delta_{f}(a)$ ($\delta_{f}>0$ lower semicontinuous).

\medskip

Consider the test function $\rho: \mathbb{C}^n \to \mathbb{R}$ given by
$$\rho(x)=\begin{cases}
      k \cdot e^{-1/(1-\|x\|^2)}, & \mbox{ if } \|x\|<1 \\
      0, & \mbox{ if } \|x\|\geq 1
    \end{cases}$$
where the constant $k$ is chosen so that $\int_{\mathbb{C}^n} \rho dm=1$. More generally, for each $\delta>0$ let $\rho_{\delta} \in \mathcal{D}(\mathbb{C}^n)$ be defined by $\rho_{\delta}(x)=\delta^{-n} \rho(x/\delta)$ for every $x \in \mathbb{C}^n$, so that $\int_{\mathbb{C}^n} \rho_{\delta} dm=1$
 and $\text{supp}(\rho_{\delta})=\bar{B}(0,\delta)$.

\medskip

Fix $\delta_0>0$, and let $U_{\delta_0}:=\{z\in U: d_U(z)> \delta_0\}$. Since $\psi$ is positive and uniformly continuous on $\overline{U_{\delta_0/2}}$, there is $\delta_1\in(0,\delta_0/2)$ such that $|\psi(a)-\psi(a-\zeta)|<\inf_{\overline{U_{\delta_0}}}\psi/2$ when $\zeta\in\bar{B}(0,\delta_1)$ and $a\in U_{\delta_0}$. Given $a\in U_{\delta_0}$, choose $b\in \mathbb{C}^n$ of small norm (namely, less than $\inf_{\overline{B}(a,\delta)} \delta_f>0$), so for $\delta\in (0,\delta_1)$,
\begin{align*}
\psi(a)/2\|b\|^2+f\ast \rho_{\delta}(a)&=\int_{\bar{B}(0,\delta)} (\psi(a)/2\|b\|^2+f(a-\zeta))\rho_{\delta}(\zeta)dm(\zeta)\\
&\leq \int_{\bar{B}(0,\delta)} (\psi(a-\zeta)\|b\|^2+f(a-\zeta))\rho_{\delta}(\zeta)dm(\zeta)\\
&\leq \int_{\bar{B}(0,\delta)} \Big(\frac{1}{2\pi}\int_0^{2\pi} f(a-\zeta+e^{i\theta}b)d\theta\Big)\rho_{\delta}(\zeta)dm(\zeta)\\
&=\frac{1}{2\pi}\int_0^{2\pi} \Big(\int_{\bar{B}(0,\delta)} f(a-\zeta+e^{i\theta}b)\rho_{\delta}(\zeta)dm(\zeta)\Big)d\theta
\end{align*}
where the last equality follows from Fubini's theorem because $f \in L^1(U, \textit{loc})$. Therefore
$$\psi(a)/2\|b\|^2+f\ast \rho_{\delta}(a)\leq \frac{1}{2\pi} \int_0^{2\pi} f\ast \rho_{\delta}(a+e^{i\theta}b)d\theta.$$
That is, $f\ast \rho_{\delta} \in C^{\infty}(U_{\delta})$ is strictly plurisubharmonic on average on $U_{\delta_0}$ for $\delta\in(0,\delta_1)$, and from the proof of Proposition \ref{pr2.1} below, we obtain
$$\sum_{j,k=1}^n \frac{\partial^2 (f\ast \rho_{\delta})(a)}{\partial z_j \partial \overline{z_k}}b_j \overline{b_k}\geq \psi(a)/2 \|b\|^2, \; \; \; \text{ for all }a\in U_{\delta_0}, b\in \mathbb{C}^n \text{ and }\delta\in(0,\delta_1).$$
Consequently, given $w\in \mathbb{C}^n$ and $\phi \in \mathcal{D}(U)$ positive, say with $\text{supp}(\phi)\subset U_{\delta_0}$, and taking $\delta_m\to 0$ with $\delta_m<\delta_1$, we have that $f\ast \rho_{\delta_m}$ converges uniformly on $\overline{U_{\delta_0}}$ to $f$, so by the dominated convergence theorem and then integration by parts, 
\begin{align*}
\int_{U_{\delta_0}} f(z) \Big(\sum_{j,k=1}^n \frac{\partial^2 \phi(z)}{\partial z_j \partial \overline{z_k}}w_j \overline{w_k}\Big)dm(z)
&=\lim_{\delta_m\to 0} \int_{U_{\delta_0}} f\ast \rho_{\delta_m}(z)\cdot \Big(\sum_{j,k=1}^n \frac{\partial^2 \phi(z)}{\partial z_j \partial \overline{z_k}}w_j \overline{w_k}\Big) dm(z)\\
&=\lim_{\delta_m\to 0} \int_{U_{\delta_0}} \Big(\sum_{j,k=1}^n \frac{\partial^2 (f\ast \rho_{\delta_m})(z)}{\partial z_j \partial \overline{z_k}}w_j \overline{w_k}\Big)\phi(z) dm(z)\\
&\geq \int_{U_{\delta_0}} \psi(z)/2 \|w\|^2 \phi(z) dm(z)
\end{align*}

i.e. $f$ is strictly plurisubharmonic in distribution.

\medskip

Now suppose that $T \in \mathcal{D}'(U)$ is a strictly plurisubharmonic distribution, where $\psi$ is a positive continuous function satisfying equation (\ref{eq2}). Then $T\ast \rho_{\delta} \in C^{\infty}(U_{\delta})$ and for all $z \in U_{\delta}$ and $b \in \mathbb{C}^n$,
\begin{align*}
\sum_{j,k=1}^n \frac{\partial^2 (T\ast \rho_{\delta})(z)}{\partial z_j \partial \overline{z_k}} b_j \overline{b_k}
&=\sum_{j,k=1}^n \frac{\partial^2 T}{\partial z_j \partial \overline{z_k}} \ast \rho_{\delta}(z) b_j \overline{b_k}\\
&=\sum_{j,k=1}^n \frac{\partial^2 T}{\partial z_j \partial \overline{z_k}}[ \rho_{\delta}(z-\cdot)] b_j \overline{b_k}\\
&\geq \int_U \psi(w)\|b\|^2 \rho_{\delta}(z-w) dm(w)\\
&=\psi\ast \rho_{\delta}(z)\|b\|^2.
\end{align*}

\medskip

Due to Proposition \ref{pr2.1}, $\psi\ast \rho_{\delta}(z)\|b\|^2/2+T\ast \rho_{\delta}(z)\leq \frac{1}{2\pi} \int_0^{2\pi} T\ast \rho_{\delta}(z+e^{i\theta}b)d\theta$ when $b$ has small norm depending lower semicontinuously on $z$ (namely, when $\|b\|< \sup\{r>0: B(z,r)\subset (\psi\ast\rho_{\delta})^{-1}(B(\psi\ast\rho_{\delta}(z),\psi\ast\rho_{\delta}(z)/2))\}$).

\medskip

Since $T$ is in particular a plurisubharmonic distribution, we get that $T\ast \rho_{\delta}$ decreases to $f \in L^1(U,\textit{loc})$ plurisubharmonic that induces $T$. Consequently, due to the dominated and monotone convergence theorems applied respectively to the positive and negative parts of each $T\ast \rho_{\delta}$, for all $z\in U$ and $b\in \mathbb{C}^n$ of small norm (namely, when $\|b\|<\sup\{r>0: B(z,r)\subset \psi^{-1}(B(\psi(z),\psi(z)/2))\}$),
\begin{align*}
\frac{1}{2\pi}\int_0^{2\pi} f(z+e^{i\theta}b)d\theta
&=\lim_{\delta_m\to 0} \frac{1}{2\pi}\int_0^{2\pi} T\ast \rho_{\delta_m}(z+e^{i\theta}b)d\theta\\
&\geq \lim_{\delta_m \to 0} (T\ast \rho_{\delta_m}(z)+\psi\ast \rho_{\delta_m}(z) \|b\|^2/2)\\
&=f(z)+\psi(z)\|b\|^2/2.
\end{align*}

The proposition is now clear.

\end{proof}

Due to Proposition \ref{thm4.1} we have that, in $\mathbb{C}^n$, a $C^2$ function is strictly plurisubharmonic if and only if it is 
strictly plurisubharmonic on average, and we will see in Proposition \ref{pr2.1} below that in infinite dimension, a $C^2$ function is strictly plurisubharmonic continuously if and only if it is 
strictly plurisubharmonic on average continuously. Meanwhile, it had been known that in finite dimension a $C^2$ function satisfies strict plurisubharmonicity if and only if in a neighborhood of each point in the domain there exists an $\epsilon>0$ such that $F-\epsilon \|\cdot\|^2$ is plurisubharmonic, but this is no longer true in infinite dimension. We constructed the following example with the help of Santill\'an Zer\'on:

\begin{example}
Let $X$ be the Hilbert space of complex sequences $(z_n)\subset \mathbb{C}^{\mathbb{N}}$ such that
$$\|(z_n)\|:=\sqrt{\sum_{k=1}^{\infty}|z_k|^2/k^2}<\infty.$$
Then the $C^2$ function on $X\setminus \{0\}$, $F(z)=\sum_{k=1}^{\infty} |z_k|^2/k^3$, has positive definite complex Hessian at $z\neq 0$,
$$D'D''F(z)(w,w)=\sum_{k=1}^{\infty}|w_k|^2/k^3, \mbox{ for all }w\in X;$$
however $F$ does not admit $\epsilon>0$ such that, for $z$ near to $z_0\neq 0$, we have plurisubharmonicity of the function
$$F(z)-\epsilon \|z\|^2=\sum_{k=1}^{\infty}|z_k|^2/k^3-\epsilon\sum_{k=1}^{\infty}|z_k|^2/k^2=\sum_{k=1}^{\infty}|z_k|^2(1/k^3-\epsilon/k^2),$$
since its complex Hessian at each $z_0\neq 0$ is given by
$$D'D''(F-\epsilon \|\cdot\|^2)(z_0)(w,w)=\sum_{k=1}^{\infty}|w_k|^2(1/k^3-\epsilon/k^2) \mbox{ for }w\in X,$$
that eventually has eigenvalues which are negative.
\end{example}

Conversely, if $X$ is infinite dimensional, a function $F:U\subset X\to[-\infty,\infty)$ that locally admits an $\epsilon>0$ such that $F-\epsilon\|\cdot\|^2$ is plurisubharmonic may not be strictly plurisubharmonic on average continuously, as exhibited by the next example.

\begin{example}
The squared norm in $\ell_{\infty}$, $\|\cdot\|^2_{\infty}$, is not strictly plurisubharmonic on average continuously since for $a=(1,0,0,\cdots)$ and $b=(0,1,0,0,\cdots)$,
$$\frac{1}{2\pi}\int_{0}^{2\pi} (\|a+e^{i\theta}b\|_{\infty}^2-\|a\|_{\infty}^2)d\theta=0,$$
however, for $\epsilon \in (0,1)$, $(1-\epsilon)\|\cdot\|_{\infty}^2$ is known to be plurisubharmonic.
\end{example}

The proof of the next result has been simplified using arguments in an old work of Takeuchi \cite{T}, yet we believe that our remaining ideas still deserve attention.

\begin{proposition}\label{pr2.1}
Let $U$ be an open domain in a Banach space $X$. A function $f \in C^2(U; \mathbb{R})$ is strictly plurisubharmonic continuously if and only if it is strictly plurisubharmonic on average continuously.
\end{proposition}
\begin{proof}
Suppose that there exists a positive function $\varphi \in C(U)$ satisfying
\begin{equation}\label{eq2.4}
\frac{1}{2\pi}\int_0^{2\pi} (f(a+e^{i\theta}b)-f(a))d\theta\geq \varphi(a)\|b\|^2.
\end{equation}
for  $a\in U$ and $b\in X$ of small norm (with size lower semicontinously depending on $a$). Given $a \in U$, fix $b \in X$ of small enough norm so that (\ref{eq2.4}) holds and we have $a+\overline{\mathbb{D}} b \subset U$, and consider the function $u(z)=f(a+z \cdot b)$, which is defined on $\mathbb{D}$.

\medskip

Then, for all $r\in (0,1)$,
\begin{equation*}
\varphi(a)\|b\|^2\cdot r^2 \leq\frac{1}{2\pi}\int_0^{2\pi} (u(r\cdot e^{i\theta})-u(0))d\theta
\end{equation*}
Consequently, by Lemma 5 in \cite{T}, and exercises 35.B and 35.D in \cite{M},
\begin{equation}\label{eq2.5}
\varphi(a)\|b\|^2 \leq \frac{\partial^2 u}{\partial z \partial \bar{z}}(0)=D'D''f(a)(b,b).
\end{equation}

\medskip

Now suppose that there exists a positive function $\varphi \in C(U)$ such that (\ref{eq2.5}) holds for all $a \in U$ and $b\in X$. Fix $a\in U$. Since $\varphi$ is continuous at $a$, there exists an upper bound $\delta(a)>0$ for the norm of $b$ to make $|\varphi(a)-\varphi(a+b)|< \varphi(a)/2$ hold ($\delta$ is lower semicontinuous for $\delta(a)=\sup\{r>0: B(a,r)\subset \varphi^{-1}(B(\varphi(a),\varphi(a)/2))\}$). Fix $b$ bounded as before, and define $M(r)=\frac{1}{2\pi}\int_0^{2\pi}[f(a+re^{i\theta}b)-f(a)]d\theta$, for all $r \in (0,1]$. Consider also the function $u(\zeta)=f(a+\zeta b)$ defined on a disk $D(0,R)\supset \overline{\mathbb{D}}$. Then, for all $\zeta \in D(0,R)$,
$$\frac{\partial^2 u}{\partial x^2}(\zeta)+\frac{\partial^2 u}{\partial y^2}(\zeta)=4\frac{\partial^2 u}{\partial \zeta \partial \bar{\zeta}}(\zeta)=4 \cdot D'D''f(a+\zeta b)(b,b)\geq 4 \cdot \varphi(a+\zeta b) \|b\|^2.$$

\medskip

Since $\frac{\partial^2 u}{\partial x^2}+\frac{\partial^2 u}{\partial y^2}=\frac{\partial^2 u}{\partial r^2}+\frac{1}{r}\frac{\partial u}{\partial r}+\frac{1}{r^2}\frac{\partial^2 u}{\partial \theta^2}$, then for $r\in(0,1)$,
$$\frac{1}{2\pi}\int_0^{2\pi}\Big(\frac{\partial^2}{\partial r^2}+\frac{1}{r}\frac{\partial}{\partial r}+\frac{1}{r^2}\frac{\partial^2}{\partial\theta^2}\Big)u(r e^{i\theta})d\theta\geq 4\cdot \frac{1}{2\pi}\int_0^{2\pi} \varphi(a+re^{i\theta} b)\|b\|^2d\theta\geq 2\cdot\varphi(a)\|b\|^2,$$
i.e. $M''(r)+\frac{1}{r}M'(r)\geq 2\cdot \varphi(a)\|b\|^2, \; \forall r \in (0,1)$.

\medskip

Thus $\Big(rM'(r)-2\cdot \varphi(a)\|b\|^2\frac{r^2}{2}\Big)'=rM''(r)+M'(r)-2\cdot \varphi(a)\|b\|^2 r\geq 0$ for all $r \in (0,1)$, so $r\Big(M'(r)-\varphi(a)\|b\|^2 r\Big)$ is an increasing function of $r$. Since clearly $r\Big(M'(r)-\varphi(a)\|b\|^2 r\Big)\to 0$ as $r\to 0$ (because $M'$ is a bounded function on $(0,\epsilon)$ for some $\epsilon>0$), we conclude that $r\Big(M'(r)-\varphi(a)\|b\|^2 r\Big)\geq 0$ for every $r \in (0,1)$. Hence $\Big(M(r)-\varphi(a)\|b\|^2 \frac{r^2}{2}\Big)' \geq 0$ for every $r>0$, so $M(r)-\varphi(a)\|b\|^2 \frac{r^2}{2}$ is an increasing function of $r$. Since clearly $M(r)-\varphi(a)\|b\|^2 \frac{r^2}{2}\to 0$ as $r\to 0$ then $M(r)\geq \varphi(a)\|b\|^2 \frac{r^2}{2}$ for each $r \in (0,1)$.

\medskip

Since $M$ is continuous on $(0,1]$, we conclude that $M(1)\geq \frac{\varphi(a)}{2}\|b\|^2 $, so indeed
$$\frac{1}{2\pi} \int_0^{2\pi} [f(a+e^{i\theta}b)-f(a)]d\theta\geq \frac{\varphi(a)}{2}\|b\|^2.$$ 

\end{proof}

An important remark about the proof of Proposition \ref{pr2.1} is that, when we have $f\in C^2(U;\mathbb{R})$ satisfying $D'D''f(a)(b,b)\geq L\|b\|^2$ for some $L>0$ we can conclude that $\frac{1}{2\pi}\int_0^{2\pi}(f(a+e^{i\theta}b)-f(a))d\theta\geq L\|b\|^2$ for all $a\in U$ and $b\in X$ with $\|b\|<d_U(a)$.

\medskip

Let us finish this section exhibiting functions strictly plurisubharmonic on average continuously, without two degrees of differentiability. For that let us discuss a family of Banach spaces $X$ having a norm $\|\cdot\|$ that is strictly plurisubharmonic on average continuously, of which some lack two degrees of differentiability \cite{DGZ}.

\medskip

The following notion of uniform convexity for complex quasi-normed spaces found in \cite{DGT} generalizes uniform \textit{c}-convexity as defined by Globevnik \cite{Gl}.

\begin{definition} If $0<q<\infty$ and $2\leq r <\infty$, a continuously quasi-normed space $(X, \| \cdot\|)$ is called \textit{$r$-uniformly PL-convex} if there exists $\lambda>0$ such that
$$\Big(\frac{1}{2\pi}\int_{0}^{2\pi} \|a+e^{i\theta}b\|^q d\theta\Big)^{1/q}\geq (\|a\|^r+\lambda\|b\|^r)^{1/r}$$
for all $a$ and $b$ in $X$; we shall denote the largest possible value of $\lambda$ by $I_{r,q}(X)$.
\end{definition}

It is known that the previous definition does not depend on $q$. Let us recall that a quasi-normed space $(X, \| \cdot\|)$ is continuously quasi-normed if $\| \cdot \|$ is uniformly continuous on the bounded sets of $X$. Banach spaces are obviously continuously quasi-normed.

\medskip

Davis, Garling and Tomczak-Jaegermann proved that for $p \in [1,2]$, $L_p(\Sigma, \Omega, \mu)$ is $2$-uniformly PL-convex (\cite[Cor.~4.2]{DGT}), and they obtain that $I_{2,p}(L_p)=I_{2,p}(\mathbb{C})\geq 1>1/2=I_{2,1}(\mathbb{C})$. Other examples of $2$-uniformly PL-convex spaces include the dual of any $C^*$-algebra (\cite[Thm.~4.3]{DGT}), the complexification of a Banach lattice with $2$-concavity constant $1$ (\cite[Cor.~4.2]{Lee}) and the non commutative $L^p(M)$, $1\leq p \leq 2$, where $M$ is a von Neumann algebra acting on a separable Hilbert space (\cite[Thm.~4]{F}.

\medskip

The following proposition gives us in particular that an $L_p(\Sigma, \Omega, \mu)$ space, for $p \in [1,2]$, has strictly plurisubharmonic norm.

\begin{proposition}\label{thm5.1}
If $(X, \| \cdot\|)$ is a 2-uniformly PL-convex Banach space then the norm $\|\cdot\|$ is strictly plurisubharmonic on average continuously.
\end{proposition}
\begin{proof}
Let $a \in X$ and $b\in B_X$. Then,
$$(\|a\|^2+I_{2,1}(X) \|b\|^2)^{1/2}\leq \frac{1}{2\pi}\int_{0}^{2\pi} \|a+e^{i\theta}b\| d\theta$$

Let $\varphi(a)=\sqrt{I_{2,1}(X)+\|a\|^2}-\|a\|>0$. Since $\|b\|<1$, 
$$\|a\|+\varphi(a) \|b\|^2 \leq (\|a\|^2+I_{2,1}(X) \|b\|^2)^{1/2}.$$ 

Thus,
$$\varphi(a) \|b\|^2+\|a\| \leq \frac{1}{2\pi}\int_{0}^{2\pi} \|a+e^{i\theta}b\| d\theta.$$
as desired.
\end{proof}

One clearly gets that the norm of a $2$-uniformly PL-convex Banach space is not only strictly plurisubharmonic on average continuously, but on bounded balls it is moreover uniformly plurisubharmonic on average:

\begin{definition} 
An upper semicontinuous function $g:U\subset X \to [-\infty,\infty)$ will be called \textit{uniformly plurisubharmonic on average} if there exists a constant $L>0$ such that for all $a\in U$ and $b\in X$ of small norm (with size depending lower semicontinuously on $a$),
$$L\|b\|^2+g(a)\leq \frac{1}{2\pi}\int_{0}^{2\pi} g(a+e^{i\theta}b)d\theta.$$

If $g:U\subset X \to [-\infty,\infty)$ is as just described above, let us specifically call it \textit{$L$-uniformly plurisubharmonic on average}.
\end{definition}

The reader can correspondingly define uniform plurisubharmonicity in distribution and prove a relationship to uniform plurisubharmonicity on average analogous to the one shown in Proposition \ref{thm4.1}.

\medskip

A last remark about functions strictly plurisubharmonic on average continuously is that if we add to one of those a plurisubharmonic function, we of course preserve a function strictly plurisubharmonic on average continuously.

\section{Strict pseudoconvexity}

Let us step back to finite dimension to discuss strict pseudoconvexity even in the case when the boundary of a given domain lacks two degrees of smoothness.

\begin{definition}
A domain $U\subset \mathbb{C}^n$ with $C^2$ boundary is called \textit{strictly pseudoconvex} when the strict Levi condition is satisfied by a $C^2$ defining function $r$ of the boundary, i.~e. when $r$ is a $C^2$ real-valued function defined on a neighborhood $V$ of $bU$ such that $U\cap V=\{z\in V: r(z)<0\}$ and $Dr(w)\neq 0$ for $w\in bU$, it holds that:
\begin{equation}\label{equation3.1}
D'D''r(w)(b,b)>0 \mbox{ when } w\in bU \mbox{ and } b\in\mathbb{C}^n \mbox{ satisfy } D'r(w)(b)=0.
\end{equation}
\end{definition}

According to \cite[Ch.~II, \S 2.8]{R}, a bounded domain $U$ in $\mathbb{C}^n$ with $C^2$ boundary is strictly pseudoconvex if and only if there is a $C^2$ defining function of $bU$, $r:V \to \mathbb{R}$, admitting a positive constant $L$ such that $D'D''r(z)(b,b)\geq L \|b\|^2$ for all $z \in V$ and $b\in\mathbb{C}^n$, where $V$ can be taken as a neighborhood of $\overline{U}$.

\medskip

Actually, it is well-known that $U$ is convex if and only if $-\log d_U$ is convex, while $U$ is pseudoconvex if and only if $-\log d_U$ is plurisubharmonic on $U$. The following is a similar result for bounded strictly pseudoconvex domains in $\mathbb{C}^n$ with $C^2$ boundary, that we proved with the help of Ramos Pe\'on.

\begin{proposition}\label{thm3.1}
Let $U$ be a bounded open domain in $\mathbb{C}^n$ with $C^2$ boundary. Then $U$ is strictly pseudoconvex if and only if there exist a positive constant L, a neighborhood $V$ of $bU$ and $\rho \in C^2(V)$ a defining function of $bU$ such that,
\begin{equation}\label{eq3.4}
D'D''(-\log|\rho|)(a)(b,b)\geq \frac{L}{|\rho(a)|} \|b\|^2 \text{ for all } a\in U\cap V \text{ and } b\in \mathbb{C}^n.
\end{equation}
Equivalently, $V$ above can be replaced by a neighborhood of $\overline{U}$.
\end{proposition}
\begin{proof}

Suppose that there exist a positive constant L, a neighborhood $V$ of $bU$ and $\rho \in C^2(V)$ a defining function of $bU$ such that $D'D''(-\log |\rho|)(a)(b,b)\geq \frac{L}{|\rho(a)|} \|b\|^2 \text{ for every } a\in U\cap V \text{ and } b\in \mathbb{C}^n$. Since for $a\in U\cap V$ and $b \in \mathbb{C}^n$ arbitrary we have
$$D'D'' (-\log |\rho|)(a)(b,b)=\frac{1}{|\rho(a)|}D'D'' \rho(a)(b,b)+\frac{1}{\rho(a)^2}\cdot |D'\rho(a)(b)|^2$$
we obtain that 
$$D'D''\rho(a)(b,b)\geq L\|b\|^2 \text{ when } a\in U\cap V \text{ and } b\in \mathbb{C}^n \mbox{ satisfy }D'\rho(a)(b)=0.$$

A passage to the limit shows that on the boundary we have what we desired:
$$D'D'' \rho(w)(b,b)>0 \text{ when } w\in bU \text{ and } b\neq 0 \text{ satisfy } D'\rho(w)(b)=0.$$

\medskip

Now suppose that $U$ is strictly pseudoconvex. Then we can find a positive constant $L$, a neighborhood $V$ of $bU$ and $\rho \in C^2(V)$ a defining function of the boundary of $U$ such that
$$D'D''\rho(a)(b,b)\geq L \|b\|^2 \text{ for all } a\in V \text{ and } b\in \mathbb{C}^n.$$

Then for $a \in U\cap V$ and $b \in \mathbb{C}^n$ arbitrary,
\begin{align*}
D'D''(-\log|\rho|)(a)(b,b)&=\frac{1}{|\rho(a)|}D'D''\rho(a)(b,b)+\frac{1}{|\rho(a)|^2}|D'\rho(a)(b)|^2\\
&\geq  \frac{L}{|\rho(a)|}\|b\|^2,
\end{align*}
as desired.
\end{proof}

Because of the previous proposition and the developments in Section 2, we consider the following concept:

\begin{definition}
Given $\ell\geq 1$ we will say that a bounded domain $U$ with $C^{\ell}$ boundary, in a Banach space $X$, is \textit{$\ell$-strictly pseudoconvex} if there exist a positive constant $L$, a neighborhood $V$ of $bU$ and $\rho\in C^{\ell}(V)$ a defining function of $bU$ such that for all $a\in U\cap V$ and $b\in X$ of small norm (with size lower semicontinuously depending on $a$),
\begin{equation}\label{equa3.2}
\frac{1}{2\pi}\int_0^{2\pi}-\log|\rho|(a+e^{i\theta}b) d\theta\geq -\log|\rho|(a)+\frac{L}{|\rho(a)|}\|b\|^2.
\end{equation}
\end{definition}

It is easy to check that an $\ell$-strictly pseudoconvex domain $U$ admits a plurisubharmonic function $\sigma$ defined on all of $U$ such that $\sigma(z)\to \infty$ as $z\to bU$ (see \cite[Ch.II, \S 2.7]{R} and \cite[Ch. VIII, \S 34]{M}).

\medskip

Clearly, a bounded domain in $\mathbb{C}^n$ with $C^2$ boundary is strictly pseudoconvex if and only if it is $2$-strictly pseudoconvex, and all $2$-strictly pseudoconvex domains are $1$-strictly pseudoconvex. Moreover, it is clear that an $\ell$-strictly pseudoconvex domain in $\mathbb{C}^n$ admits an associated function $\sigma$ as a plurisubharmonic exhaustion function, hence $\ell$-strictly pseudoconvex domains in $\mathbb{C}^n$ are pseudoconvex. Likewise, an $\ell$-strictly pseudoconvex domain in a Banach space $X$ is pseudoconvex because its restriction to each finite-dimensional subspace admits a restricted associated function $\sigma$ as a plurisubharmonic exhaustion function.

\medskip

Let us now look at the following necessary condition for $\ell$-strict pseudoconvexity in $\mathbb{C}^n$.

\begin{proposition}\label{pr3.2}
If a domain $U\subset\mathbb{C}^n$ is $\ell$-strictly pseudoconvex for some $\ell \geq 1$ then there exist a positive constant $L$, a neighborhood $V$ of $bU$ and a defining function $\rho\in C^{\ell}(V)$ of $bU$ such that for all $a\in U\cap V$ and $b\in \mathbb{C}^n$ with $D'\rho(a)(b)=0$ and $\|b\|<d_{U\cap V}(a)$we have,
$$\frac{1}{2\pi}\int_0^{2\pi} \rho(a+e^{i\theta}b)\geq \rho(a)+L\|b\|^2.$$
\end{proposition}
\begin{proof}
Since $U$ is $\ell$-strictly pseudoconvex, there exist a positive constant $L_0$, a neighborhood $V$ of $bU$ and a defining function $\rho\in C^{\ell}(V)$ of $bU$ such that for all $a\in U\cap V$ and $b\in \mathbb{C}^n$ of small norm (with size lower semicontinuously depending on $a$),
$$\frac{1}{2\pi}\int_0^{2\pi}-\log|\rho|(a+e^{i\theta}b) d\theta\geq -\log|\rho|(a)+\frac{L_0}{|\rho(a)|}\|b\|^2.$$

For each $n\in\mathbb{N}$ define the function $\rho_{1/n}$ as in Proposition \ref{thm4.1}, and take $\sigma_n=(-\log(-\rho))\ast \rho_{1/n}$, which is defined and smooth on $(U\cap V)_{1/n}$. Then, as in the proof of Proposition \ref{thm4.1}, there is $n_1\geq n$ such that for all $m\geq n_1$ we have that $\sigma_m$ is strictly plurisubharmonic on average on $(U\cap V)_{1/n}$ with respect to the function $\frac{L_0}{2|\rho|}$, and because of the proof of proposition \ref{pr2.1} we have that if $L=L_0/2$, for all $a\in (U\cap V)_{1/n}$, $b\in\mathbb{C}^n$ and $m\geq n_1$,
$$D'D''\sigma_m(a)(b,b)\geq \frac{L}{|\rho|(a)}\|b\|^2.$$
Because of distribution theory it is clear that $\sigma_m$ converges uniformly on each $\overline{(U\cap V)_{\frac{2}{n}}}$ to $-\log(-\rho)$, and likewise the vector-valued function $D'\sigma_m$ converges uniformly on each $\overline{(U\cap V)_{\frac{2}{n}}}$ to $D'(-\log(-\rho))=\frac{-1}{\rho}D'\rho$. Define $r_m=-e^{-\sigma_m}$ on $(U\cap V)_{1/m}$, which again converges uniformly on each $\overline{(U\cap V)_{\frac{2}{n}}}$ to $\rho$.

\medskip

Pick $n\in\mathbb{N}$. Now choose $\epsilon_n>0$ such that $\epsilon_n<\min\Big(\frac{L}{M_n+1}, 1\Big)$, where $M_n=\sup_{a\in \overline{(U\cap V)_{2/n}}}(L/|\rho(a)|+|\rho(a)|)$.
Without loss of generality, for all $m\geq n_1$ we have that $|e^{-\sigma_m}-|\rho| \; |\leq \epsilon_n$ on $\overline{(U\cap V)_{\frac{2}{n}}}$, and that $|\; |D'\sigma_m(a)(b)|^2-|\frac{1}{|\rho(a)|}D'\rho(a)(b)|^2\; |\leq \epsilon_n\|b\|^2$ when $a\in \overline{(U\cap V)_{\frac{2}{n}}}$ and $b\in\mathbb{C}^n$.

\medskip

Then, when $a\in \overline{(U\cap V)_{2/n}}$ and $D'\rho(a)(b)=0$, we have that for $m\geq n_1$,
\begin{align*}
D'D''r_m(a)(b,b)&=e^{-\sigma_m(a)}D'D''\sigma_m(a)(b,b)-e^{-\sigma_m(a)}|D'\sigma_m(a)(b)|^2\\
&\geq e^{-\sigma_m(a)} \frac{L}{|\rho(a)|}\|b\|^2-e^{-\sigma_m(a)}|D'\sigma_m(a)(b)|^2\\
&\geq \Big(\frac{|\rho(a)|-\epsilon_n}{|\rho(a)|}\Big)L\|b\|^2-(|\rho(a)|+\epsilon_n)\Big(\Big|\frac{1}{\rho(a)}D'\rho(a)(b)\Big|^2+\epsilon_n\|b\|^2\Big)\\
&=(L-\epsilon_n\Big(\frac{L}{|\rho(a)|}+|\rho(a)|\Big)-\epsilon_n^2)\|b\|^2
\end{align*}
where $L-\epsilon_n\Big(\frac{L}{|\rho(a)|}+|\rho(a)|\Big)-\epsilon_n^2\geq L-\epsilon_n(M_n+1)>0$.
 Then, following the proof of Proposition \ref{pr2.1}, each $r_m$ with $m\geq n_1$ satisfies for $a\in (U\cap V)_{2/n}$ as well as $b$ with $D'\rho(a)(b)=0$ and $\|b\|<d_{(U\cap V)_{2/n}}(a)$ ,
$$\frac{1}{2\pi}\int_0^{2\pi} (r_m(a+e^{i\theta}b)-r_m(a))d\theta\geq (L-\epsilon_n(M_n+1))\|b\|^2.$$
so taking the limit as $m\to\infty$,
$$\frac{1}{2\pi}\int_0^{2\pi} (\rho(a+e^{i\theta}b)-\rho(a))d\theta\geq (L-\epsilon_n(M_n+1))\|b\|^2.$$
As $\epsilon_n\to 0$ and then $n\to \infty$, we conclude what we desired.
\end{proof}

Proposition \ref{pr3.2} is significant since it confirms a very expected behavior of a defining function of an $\ell$-strictly pseudoconvex domain. However, if $\ell=1$, it does not recover a condition at the boundary resembling equation (\ref{equation3.1}), so we consider it separately.

\begin{definition}
Given a domain $U\subset X$ with $C^{1}$ boundary, we will say that $U$ is \textit{$1$-strictly pseudoconvex at the boundary} if it is pseudoconvex and there exist a neighborhood $V$ of $bU$, $\rho\in C^{1}(V)$ a defining function of $bU$ and $\varphi\in C(bU)$ a positive function, such that for all $w\in bU$ and $b\in\mathbb{C}^n$ with $D'\rho(w)(b)=0$ and $\|b\|$ small (lower semicontinuously depending on $w$) we have
\begin{equation}\label{equ3.4}
\frac{1}{2\pi}\int_0^{2\pi}\rho(w+e^{i\theta}b)d\theta\geq \rho(w)+\varphi(w)\|b\|^2.
\end{equation}
\end{definition}

We will soon explore examples of domains $1$-strictly pseudoconvex at the boundary, both in finite and infinite dimension. But first let us see that the following is a sufficient condition for $\ell$-strict pseudoconvexity when $\ell\geq 1$.

\begin{proposition}
Let $\ell\geq 1$. If $U$ in a Banach space $X$ is a bounded open domain that admits a $C^{\ell}$ defining function of $bU$, $r:V\supset bU\to\mathbb{R}$, which is uniformly plurisubharmonic on average on $U\cap V$ then $U$ is $\ell$-strictly pseudoconvex.
\end{proposition}
\begin{proof}

Let $r:V\supset bU \to \mathbb{R}$ be a $C^{\ell}$ defining function of $bU$ which is uniformly plurisubharmonic on average on $U\cap V$. Then we can find $L>0$ such that, for all $a\in U\cap V$ and $b\in\mathbb{C}^n$ of small norm (with size depending lower semicontinuously on $a$), we have that $a+\overline{\mathbb{D}}b$ is contained in $U\cap V$ and
\begin{equation}\label{equation3.2}
\frac{1}{2\pi}\int_0^{2\pi} r(a+e^{i\theta}b)d\theta\geq r(a)+L\|b\|^2.
\end{equation}

\medskip

Fix $a\in U\cap V$ and $b\in X$ of small norm as before. Since $a+\overline{\mathbb{D}}\cdot b\subset U\cap V$, we have that $r(a+e^{i\theta}b)<0$ for each $\theta\in[0,2\pi]$.

\medskip

Since $x\mapsto -\log(-x)$ is convex on $(-\infty,0)$, we obtain by Jensen's inequality that
\begin{equation}\label{equation3.3}
\frac{1}{2\pi}\int_0^{2\pi}-\log(-r)(a+e^{i\theta}b)d\theta\geq -\log(-\frac{1}{2\pi}\int_0^{2\pi}r(a+e^{i\theta}b)d\theta).
\end{equation}

Meanwhile, due to equation (\ref{equation3.2}) we have that
$$0<\frac{-1}{2\pi}\int_{0}^{2\pi} r(a+e^{i\theta}b)d\theta\leq -(r(a)+L\|b\|^2),$$
so using that $-\log$ is decreasing on $(0,\infty)$ we obtain that
\begin{equation}\label{equation3.4}
-\log(\frac{-1}{2\pi}\int_{0}^{2\pi} r(a+e^{i\theta}b)d\theta)\geq -\log(-(r(a)+L\|b\|^2))
\end{equation}

As a consequence of equations (\ref{equation3.3}) and (\ref{equation3.4}) we obtain that

\begin{align*}
\frac{1}{2\pi}\int_0^{2\pi}-\log(-r)(a+e^{i\theta}b)d\theta& \geq -\log(-(r(a)+L\|b\|^2)) \\
& =-\log (-r(a)(1-\frac{L\|b\|^2}{[-r(a)]}))\\
& =-\log(-r(a)) -\log(1-\frac{L\|b\|^2}{[-r(a)]})\\
& \geq -\log(-r(a))+L\|b\|^2/[-r(a)]
\end{align*}
as desired.
\end{proof}

Note that if $U$ bounded admits a $C^1$ defining function $\rho:V\supset bU\to\mathbb{R}$ which is uniformly plurisubharmonic on average on $U\cap V$ then there exists $\delta>0$ such that for all $c\in(-\delta,0)$, $U_{\rho,c}:=(U\setminus V)\cup\{z\in V: \rho(z)< c\}$ is 1-strictly pseudoconvex, and moreover 1-strictly pseudoconvex at the boundary. 

\medskip

We proceed to present examples of domains $1$-strictly pseudoconvex at the boundary. The following concept of strong $\mathbb{C}$-linear convexity is succinctly studied for the ambient space $\mathbb{C}^n$ in \cite{LS}.

\begin{definition}
A domain $U\subset X$ is called \textit{strongly $\mathbb{C}$-linearly convex} if $bU$ is of class $C^1$ and there is a $C^1$ defining function $\rho:V\supset bU\to\mathbb{R}$ which satisfies that for a certain $C> 0$,
$$|D'\rho(w)(w-z)|\geq C\|w-z\|^2 \text{ for all }w\in bU, \; z\in \overline{U}$$
(in $\mathbb{C}^n$, if this condition holds for a defining function then it passes to every defining function).
\end{definition}

\begin{definition}
A domain $U$ in $X$ is called \textit{weakly linearly convex} if for every $w\in bU$ there exists $\pi$ affine complex hyperplane such that $w\in \pi\subset X\setminus U$.
\end{definition}

\begin{example}
Domains $U$ in a Banach space $X$ that are strongly $\mathbb{C}$-linearly convex are weakly linearly convex, since given $w\in bU$ we have that $w\in w+T_w^{\mathbb{C}}(bU)$, and for all $z\in U$, there exists $C>0$ such that
$$\inf_{b\in T_w^{\mathbb{C}}(bU)}\|w+b-z\|\geq \inf_{b\in T_w^{\mathbb{C}}(bU)} \frac{1}{\|D'\rho(w)\|} \cdot |D'\rho(w)(w+b-z)| \geq \frac{C}{\|D'\rho(w)\|} \|w-z\|^2>0.$$

It is well known that domains weakly linearly convex are holomorphically convex (see \cite[Ch. IV, \S 4.6]{H2}, whose proof still works in infinite dimension), thus domains weakly linearly convex are pseudoconvex \cite[Ch.~VIII, \S 37]{M}.

\medskip

Let us now show that domains strongly $\mathbb{C}$-linearly convex $U\subset\mathbb{C}^n$ with $C^{1,1}$ boundary admit a $C^1$ defining function $\rho: V\supset bU\to\mathbb{R}$ satisfying equation (\ref{equ3.4}) for all $w\in bU$ and $b\in\mathbb{C}^n$ with $D'\rho(w)(b)=0$ and $\|b\|<\min(d_{V}(w),1)$.

\medskip

Consider the function $$\rho(z)=\begin{cases} -d_U(z), &\text{ if }z\in U\\d_U(z), &\text{ if }z\in U^c \end{cases}.$$ Since $bU$ is $C^{1,1}$, $\rho$ is a $C^1$ defining function in a neighborhood $V$ of $bU$.

\medskip

We now fix $w\in bU$ and $b\in T_w^{\mathbb{C}}(bU)$ with $\|b\|<\min(d_{V}(w),1)$. For every $\theta\in[0,2\pi]$, we know that $d_U(w+e^{i\theta}b)=\inf_{w'\in bU}\|w+e^{i\theta}b-w'\|$, and for each $w'\in bU$ let us consider two cases:

\begin{itemize}
\item If $\|w-w'\|\leq \|b\|/2$ then $\|w+e^{i\theta}b-w'\|\geq \|b\|-\|w-w'\|\geq \|b\|/2\geq \|b\|^2/2$
\item If $\|w-w'\|\geq \|b\|/2$ then as we did before we obtain that 
$$\|w+e^{i\theta}b-w'\|\geq \inf_{b'\in T_w^{\mathbb{C}}(bU)}\|w+b'-w'\|\geq \frac{C}{\|D'\rho(w)\|} \|w-w'\|^2\geq \frac{C}{\|D'\rho(w)\|} \|b\|^2/4$$
\end{itemize}
then taking $\varphi(w)=\min\Big(\frac{1}{2}, \frac{C}{4\|D'\rho(w)\|}\Big)$, we obtain that $d_U(w+e^{i\theta}b)\geq \varphi(w) \|b\|^2$. 
Since $w+e^{i\theta}b\notin U$ for every $\theta\in[0,2\pi]$, $\rho(w+e^{i\theta}b)-\rho(w)=d_U(w+ e^{i\theta}b)\geq \varphi(w)\|b\|^2$, hence
$$\frac{1}{2\pi}\int_0^{2\pi} (\rho(w+e^{i\theta}b)-\rho(w))d\theta\geq \varphi(w)\|b\|^2.$$
\end{example}

We will come back to domains $1$-strictly pseudoconvex at the boundary, but for now let us discuss domains $U$ that are $\ell$-strictly pseudoconvex at the boundary as well as $\ell$-strictly pseudoconvex. Obviously, this will happen when a defining function of $bU$ is uniformly plurisubharmonic on average on all of its domain. Thus we will be interested in the following notion.
\begin{definition}
If $\ell\geq 1$ and $U$ in a Banach space $X$ is a bounded open domain that admits a $C^{\ell}$ defining function of $bU$ which is uniformly plurisubharmonic on average then we will say that $U$ is \textit{$\ell$-uniformly pseudoconvex}. 
\end{definition}

The reader may recall that a bounded domain in $\mathbb{C}^n$ with $C^2$ boundary is strictly pseudoconvex if and only if it is $2$-uniformly pseudoconvex. 

\medskip

Let us see nontrivial examples of $1$-uniformly pseudoconvex domains, of which at least $B_{\ell_2}$ is $2$-uniformly pseudoconvex.

\begin{proposition}\label{pr3.3}
If $n\in\mathbb{N}$ and $p \in (1,2)$, $B_{\ell_p^n}$ is $1$-uniformly pseudoconvex.
\end{proposition}
\begin{proof}

Let $n\in\mathbb{N}$ and $p \in (1,2)$. We have seen that the norm $\|\cdot\|_p$ of $\ell_p^n$ is uniformly plurisubharmonic on average, consequently so is also the function $r=\|\cdot\|_p-1$. Moreover, as shown in \cite{DGZ}, the norm $\|\cdot\|_p$ is of class $C^1$ except at $0$, and hence so is $r$. To complete showing that $bB_{\ell_p^n}$ has $1$ degree of differentiability, observe that $r$ satisfies that 
$$\{z\in (2B_{\ell_p^n})\setminus\{0\}: r(z)<0\}=B_{\ell_p^n}\setminus\{0\}$$
and the gradient of $r$ is not null at every element of $bB_{\ell_p^n}$, since for each point $z$ in $bB_{\ell_p^n}$ there is $i\in\{1,\cdots, n\}$ so that $z_i\neq 0$, and for such $i$ we have,
$$\frac{\partial |z_i|^p}{\partial z_i}\Big|_z=\frac{\partial (z_i\overline{z_i})^{p/2}}{\partial z_i}\Big|_z=p/2(z_i\overline{z_i})^{p/2-1}\overline{z_i}=p/2|z_i|^{p-2}\overline{z_i}\neq 0$$
so $\nabla (\|\cdot\|_p-1)|_z\neq 0$.
\end{proof}

\begin{proposition}\label{pr3.6}
If $p\in(1,2]$ then $B_{\ell_p}$ is $1$-uniformly pseudoconvex.
\end{proposition}
\begin{proof}
Let $p\in(1,2]$. As in the previous proposition, we know that the norm $\|\cdot\|_p$ of $\ell_p$ is uniformly plurisubharmonic on average, and hence so  is $r=\|\cdot\|_p-1$. Again, the norm $\|\cdot\|_p$ is of class $C^1$ except at $0$, and thus $r$ is $C^1$ except at $0$ as well. And $r$ satisfies that 
$$\{z\in (2B_{\ell_p})\setminus\{0\}: r(z)<0\}=B_{\ell_p}\setminus\{0\}$$
where the derivative of $r$ is not null at points $z$ in $bB_{\ell_p}$, because for each $z$ in $bB_{\ell_p}$ there exists $n\in\mathbb{N}$ such that $z_n\neq 0$, so if $e_n$ denotes the $n$-th element of the canonical basis of $\ell_p$ and $\theta=\arg(z_n) \in [0, 2\pi)$, we have that
\begin{align*}
Dr(z)(e^{i\theta}e_n)&=\lim_{h\to 0^+}\frac{\|z+he^{i\theta}e_n\|-\|z\|}{h}\\
&=\lim_{h\to 0^+}\frac{(1-|z_n|^p+(|z_n|+h)^p)^{1/p}-1}{h}\\
&=\lim_{h\to 0^+}\frac{|z_n|^{p-1}h+o(|h|)}{h},
\end{align*}
where the last equality has been obtained with Taylor series, leading us to obtain that $\|Dr(z)\|\geq |Dr(z)(e^{i\theta}e_n)|=|z_n|^{p-1}> 0$. We conclude that $B_{\ell_p}$ is $1$-uniformly pseudoconvex.\end{proof}

Suitable modifications to the defining functions $\|\cdot\|_p-1$ by plurisubharmonic functions can lead to examples of domains that are at least $1$-strictly pseudoconvex at the boundary, as in the next examples.

\begin{proposition}
If $p\in(1,2]$, then $\{z\in\mathbb{C}^n: -\text{Re}(z_1)+\|z\|_p<1 \}$ is $1$-strictly pseudoconvex at the boundary.
\end{proposition}
\begin{proof}
From Proposition \ref{thm5.1}, the function $r(z)=-\text{Re}(z_1)+\|z\|_p-1$ is strictly plurisubharmonic on average. Also, since the function $r$ is in particular a continuous plurisubharmonic function, $\{z\in\mathbb{C}^n: r(z)<0\}$ is pseudoconvex \cite[Ch.~VIII, \S 38]{M}. 

\medskip

Finally, the gradient of $r$ is not null at each $z\in b\{z'\in\mathbb{C}^n: -\text{Re}(z'_1)+\|z'\|_p<1 \}$ because, if $z_j\neq 0$ for some $j>1$ then
$$\frac{\partial r}{\partial z_j}(z)=\frac{1}{2}\|z\|_p^{1-p}|z_j|^{p-2}\overline{z_j}\neq 0,$$
or if $\text{Im}(z_1)\neq 0$ then finding Taylor series one gets that
$$\frac{\partial r}{\partial y_1}(z)=\|z\|_p^{1-p}|z_1|^{p-2}\text{Im}(z_1)\neq 0;$$
otherwise $\text{Im}(z_1)=z_2=\cdots=z_n=0$, and since $-\text{Re}(z_1)+\|z\|_p=1$, we have that $\text{Re}(z_1)=-1/2$, where
$$\frac{\partial r}{\partial x_1}(-\frac{1}{2} e_1)=-1+\|z\|_p^{1-p}|z_1|^{p-2}\text{Re}(z_1)\big|_{(-\frac{1}{2}e_1)}=-1+(\frac{1}{2})^{1-p}(\frac{1}{2})^{p-2}(-\frac{1}{2})=-2.$$
\end{proof}

\begin{proposition}
If $p\in(1,2]$, then $\{z\in\ell_p: -\text{Re}(z_1)+\|z\|_p<1 \}$ is $1$-strictly pseudoconvex at the boundary.
\end{proposition}

The proof is analogous to the finite-dimensional one.

\medskip

Due to distribution theory, given $\ell\geq 1$, a bounded open domain $U\subset\mathbb{C}^n$ is $\ell$-uniformly pseudoconvex if and only if there exist open domains $V_1$ and $V_2$ such that $V_2\supset \overline{V_1} \supset bU$, there exist positive constants $L$, $M$ and $M'$ and a sequence of pointwise bounded $C^{\ell+1}$ functions on $V_2$, $\{r_m\}$, which are $L$-uniformly plurisubharmonic on average on $V_2$ and such that $M'\geq \|\nabla r_m(z)\|\geq M$ for all $z\in \overline{V_1}$, and there exists a $C^{\ell}$ function, $r:V_2\to\mathbb{R}$, satisfying $U\cap V_2=\{z\in V_2:r(z)<0\}$ and whenever $\alpha$ is a multi-index with $0\leq |\alpha|\leq \ell$ we have that the following holds uniformly on $\overline{V_1}\subset \mathbb{R}^{2n}$:
$$\frac{\partial^{\alpha}r}{\partial x^{\alpha}}=\lim_{m\to\infty} \frac{\partial^{\alpha}r_m}{\partial x^{\alpha}}.$$
Equivalently, $\{r_m\}$ above can be replaced by a family of $C^{\infty}$ functions.

\medskip

Based on the previous remark and the Arzel\`{a}-Ascoli theorem we define what comes next.

\begin{definition}
Let us call $U\Subset\mathbb{C}^n$ a \textit{$0$-uniformly pseudoconvex} domain when there exist positive constants $L$, $M$ and $M'$, bounded neighborhoods $V_1$ and $V_2$ of $bU$ with $V_2\supset \overline{V_1}$ and a function $r:V_2\to\mathbb{R}$ such that $U\cap V_2=\{z\in V_2: r(z)<0\}$ and $r$ is the limit of a sequence of pointwise bounded $C^1$ functions given on $V_2$, $\{r_m\}$, such that each $r_m$ is $L$-uniformly plurisubharmonic on average on $V_2$ and $M'\geq \|\nabla r_m(z)\|\geq M$ for all $z\in \overline{V_1}$.
\end{definition}

As before, it is not hard to check that $0$-uniformly pseudoconvex domains are pseudoconvex, and obviously $\ell$-uniformly pseudoconvex domains are $0$-uniformly pseudoconvex. 

\medskip

Even though in finite dimension we have that strong and strict pseudoconvexity are the same notion, for infinite dimension we now define strong pseudoconvexity and leave open the general concept of strict pseudoconvexity.

\begin{definition}
We will say that a bounded open set $U$ in a Banach space $X$ is \textit{strongly pseudoconvex} if $U\cap M$ is $0$-uniformly pseudoconvex for each finite-dimensional subspace $M$ of $X$. 
\end{definition}

\begin{theorem}
$B_{\ell_1}$ is strongly pseudoconvex.
\end{theorem}
\begin{proof}
We have that $r(z)=\|z\|_1-1$ on $2B_{\ell_1}\setminus\{0\}$ satisfies that $B_{\ell_1}\setminus\{0\}=\{z \in 2B_{\ell_1}\setminus\{0\}: r(z)<0\}$. Also, $r$ is the pointwise limit of $\| \cdot\|_{p_n}-1$ for $2> p_n\to 1^+$, where each $\| \cdot\|_{p_n}-1$ is $(\sqrt{\frac{1}{e}+4}-2)$-uniformly plurisubharmonic on average, because for $p\in(1,2)$, $I_{2,p}(L_p)\geq 1$ and we can argue as in the proofs of \cite[Thm.~2.4]{DGT} and Proposition \ref{thm5.1}, indeed, for $z\in 2B_{\ell_1}\setminus\{0\}$ and $b\in B_{\ell_1}$,
\begin{align*}
\frac{1}{2\pi}\int_0^{2\pi} \|z+e^{i\theta}b\|_p d\theta &\geq (\frac{1}{2\pi}\int_0^{2\pi} \|z+\frac{1}{\sqrt{e}}e^{i\theta}b\|_p^2)^{1/2}\\
&\geq (\frac{1}{2\pi}\int_0^{2\pi} \|z+\frac{1}{\sqrt{e}}e^{i\theta}b\|_p^p)^{1/p}\\
&\geq (\|z\|_p^2+\frac{1}{e}\|b\|_p^2)^{1/2}\\
&\geq \|z\|_p+(\sqrt{\frac{1}{e}+4}-2)\|b\|_p^2.
\end{align*}
And the $C^1$ functions $\{|\|\cdot\|_p-1|\}_{p\in(1,2)}$ on $2B_{\ell_1}\setminus\{0\}$ are bounded by $1$ because $\|\cdot\|_p\leq\|\cdot\|_1$. Clearly the same holds when we restrict those functions to finite-dimensional subspaces.

\bigskip

Now let $M$ be a finite-dimensional subspace of $\ell_1$. Then, since all norms are equivalent in $M$, we will just bound the operator norm of $D(\|\cdot\|_p-1)(z)$ for $p\in(1,2)$ and $z\in \overline{2B_M}\setminus\{0\}$, where $\|z\|_1\leq C_M \|z\|_2$ for all $z\in M$, for a constant $C_M\geq1$ only depending on $M$, so
$$\|D (\|\cdot\|_p-1)(z)\|\geq \frac{D (\|\cdot\|_p-1)( z)(z)}{\|z\|_1}=\lim_{h\to 0^+}\frac{\|z+hz/\|z\|_1\|_p-\|z\|_p}{h}=\|z\|_p/\|z\|_1$$
which is above $\|z\|_2/\|z\|_1\geq 1/C_M$; and to bound the gradient from above, 
$$\|D (\|\cdot\|_p-1)\|=\sup_{w\neq 0} \frac{|D (\|\cdot\|_p-1)(z)(w)|}{\|w\|_1}$$
and since for $e_n$ the $n$-th element of the canonical basis of $\ell_1$ and $\theta_n=\text{arg}(z_n)$ we have, as in the proof of Proposition \ref{pr3.6} that 
\begin{align*}
D(\|\cdot\|_p-1)(z)(e^{i\theta_n}e_n)&=\|z\|_p^{1-p}|z_n|^{p-1},\\
D(\|\cdot\|_p-1)(z)(ie^{i\theta_n}e_n)&=0,
\end{align*}
then, if $<\cdot,\cdot>$ represents the usual inner product in $\mathbb{R}^2$, for any $w\neq 0$ in $M$,
\begin{align*}
|D (\|\cdot\|_p-1)(z)(w)|&=\Big|\sum_n <e^{i\theta_n}, w_n> D(\|\cdot\|_p-1)(z)(e^{i\theta_n}e_n)\Big|\\
&\leq \sum_n |w_n| \|z\|_p^{1-p}|z_n|^{p-1}\\
&\leq \|z\|_p^{1-p} \|z\|_{\infty}^{p-1} \|w\|_1\\
&\leq \|z\|_2^{1-p} \|z\|_1^{p-1}\|w\|_1
\end{align*}
hence
$$\|D (\|\cdot\|_p-1)(z)\|_2\leq \|z\|_2^{1-p} \|z\|_1^{p-1}\leq (C_M)^{p-1}\leq C_M.$$
\end{proof}

\begin{corollary}\label{thm3.4}
$B_{\ell_1^n}=\{z\in\mathbb{C}^n:\sum_{j=1}^n |z_j|<1\}$ is $0$-uniformly pseudoconvex.
\end{corollary}

We continue our discussion of $\ell$-uniform pseudoconvexity with a characterization of $0$-uniformly pseudoconvex domains as a certain limit of $1$-uniformly pseudoconvex domains.

\begin{proposition}\label{pr3.1}
In $\mathbb{C}^n$, a domain $U$ is $0$-uniformly pseudoconvex if and only if it is exhausted by an increasing sequence $\{U_m\}$ of $1$-uniformly pseudoconvex domains given by a respective family of $C^1$ defining functions  $r_m: V'\supset \overline{V}\supset bU_m \cup bU \to \mathbb{R}$ such that $\infty>\limsup |r_m(z)|\geq \liminf |r_m(z)|>0$ for each $z\in \overline{V}\setminus bU$, and such that there exist common positive bounds $L$, $M$ and $M'$ satisfying that each $r_m$ is $L$-uniformly plurisubharmonic on average on $V$ and $M' \geq \|Dr_m(z)\|\geq M$ for all $z \in \overline{V}$ compact and $w\in\mathbb{C}^n$.
\end{proposition}
\begin{proof}
Suppose that we can exhaust $U$ with an increasing sequence $\{U_m\}$ of $1$-uniformly pseudoconvex domains, where each $bU_m$ is given by the $C^1$ defining function $r_m:V'\supset\overline{V}\supset bU_m\cup bU\to\mathbb{R}$ so that the family $\{r_m\}$ satisfies $\infty>\limsup |r_m(z)|\geq \liminf |r_m(z)|>0$ for each $z\in \overline{V}\setminus bU$, each $r_m$ is $L$-uniformly plurisubharmonic on average on $V$, and for all $z$ in the compact set $\overline{V}$ as well as any $w\in\mathbb{C}^n$, it holds that $M' \geq \|Dr_m(z)\|\geq M$. Then we can use Arzel\`a-Ascoli theorem to find a subsequence $\{r_{m_k}\}$ of $\{r_m\}$ as restricted to $\overline{V}$ that converges uniformly to a function $r:\overline{V}\to\mathbb{R}$, and it will also hold that $U\cap V=\{z\in V: r(z)<0\}$.
\medskip

Conversely, suppose that $r$ is the uniform limit on $\overline{V}$ compact of a sequence of $C^1$ functions on $V'\supset \overline{V}$, $\{r_m\}$, that are $L$-uniformly plurisubharmonic on average and such that $M'\geq \|Dr_m(z)\|\geq M$ for all $z\in \overline{V}$, and where $r: V'\supset \overline{V}\supset bU \to\mathbb{R}$ satisfies $U\cap V'=\{z\in V':r(z)<0\}$. Let $V_0$ be a neighborhood of $bU$ with $\overline{V_0}\subset V$. For each $m\in\mathbb{N}$ take $c_m=\min_{\overline{V}\setminus U} r_m$, and set $U_m=(U\setminus \overline{V_0})\cup\{z\in V: r_m(z)-c_m<0\}$. It is clear that $c_m \to 0$ and that $\cup U_m=U$. By subtracting small positive numbers $d_m$ to each $c_m$  with $d_m\to 0$ and then passing to a subsequence, we may assume that $\{U_m\}$ is an increasing sequence of open sets whose union is $U$, where $bU_m$ has $C^1$ defining function $r_m-(c_m-d_m)$ defined on $V_0$. We conclude that $U$ satisfies the desired conditions. 
\end{proof}

For $\ell\geq 1$, the reader can similarly show a characterization of $\ell$-uniformly pseudoconvex domains in $\mathbb{C}^n$ as a certain limit of $(\ell+1)$-uniformly pseudoconvex domains. In particular, the following holds.

\begin{proposition}
A domain $U$ in $\mathbb{C}^n$ is $1$-uniformly pseudoconvex if and only if it is exhausted by an increasing sequence $\{U_m\}$ of strongly pseudoconvex domains given by a respective family of $C^2$ defining functions $r_m: V'\supset \overline{V}\supset bU_m \cup bU \to \mathbb{R}$ such that $\infty>\limsup |r_m(z)|\geq \liminf |r_m(z)|>0$ for each $z\in \overline{V}\setminus bU$, that each family $\{\frac{\partial r_m}{\partial x_i}\}$ is equicontinuous on $\overline{V}$ for $i=1, \cdots, 2n$, and that there exist common positive bounds $L$, $M$ and $M'$ such that $D'D''r_m(z)(w,w)\geq L \|w\|^2$ and $M' \geq \|Dr_m(z)\|\geq M$ for all $z \in \overline{V}$ compact and $w\in\mathbb{C}^n$.\\
Equivalently, $\{r_m\}$ above can be replaced by a family of $C^{\infty}$ functions.
\end{proposition}

Let us use the previous approximation result to show an invariance property of $1$-uniformly pseudoconvex domains in $\mathbb{C}^n$, which in particular yields that the image of a $1$-uniformly pseudoconvex domain in $\mathbb{C}^n$ under an affine isomorphism remains $1$-uniformly pseudoconvex. Of course, this extends our family of examples of $1$-uniformly pseudoconvex domains.

\begin{proposition}
The biholomorphic image of a $1$-uniformly pseudoconvex domain $U$ in $\mathbb{C}^n$ is still $1$-uniformly pseudoconvex, when the domain of the biholomorphism contains $\overline{U}$.
\end{proposition}
\begin{proof}
Suppose that $U$ is a $1$-uniformly pseudoconvex domain which is exhausted by the increasing sequence $\{U_m\}$ of strongly pseudoconvex domains given by the respective family of $C^2$ defining functions $r_m: V'\supset \overline{V}\supset bU_m \cup bU \to \mathbb{R}$ such that $\infty>\limsup |r_m(z)|\geq \liminf |r_m(z)|>0$ for each $z\in \overline{V}\setminus bU$, for which there exist common positive bounds $L$, $M$ and $M'$ such that $D'D''r_m(z)(w,w)\geq L \|w\|^2$ and $M' \geq \|Dr_m(z)\|\geq M$ for all $z$ in the compact set $\overline{V}$ and $w\in\mathbb{C}^n$, and such that each family $\{\frac{\partial r_m}{\partial x_i}\}$ is equicontinuous, $i=1, \cdots, 2n$.

\medskip

Then, as in \cite[Ch.~II, \S 2.6]{R}, if $F:W'\supset \overline{W}\supset U \to\mathbb{C}^n$ is a biholomorphic map, then $F(U)$ is exhausted by the increasing sequence of domains $\{F(U_m)\}$ given by the respective family of defining functions $\rho_m=r_m\circ F^{-1}: F(W'\cap V')\supset F(\overline{W}\cap\overline{V})\supset bF(U_m)\cup bF(U)\to\mathbb{R}$, which satisfy for $z\in\overline{W}\cap\overline{V}$ and $w\in\mathbb{C}^n$,
\begin{align*}
i)&\infty>\limsup |\rho_m\circ F(z)|\geq \liminf |\rho_m\circ F(z)|>0\text{ when }z\notin bU,\\
ii)&D'D''r_m(z)(w,w)=D'D''\rho_m(F(z))(F'(z)w, F'(z)w),\\
iii)&Dr_m(z)(w)=D\rho_m(F(z))(F'(z)w);
\end{align*}
where each $F'(z)$ is a nonsingular $\mathbb{C}$-linear map, so it defines an isomorphism on $\mathbb{C}^n$. Hence each domain $F(U_m)$ is still strongly pseudoconvex with
 \begin{align*}
 i')&D'D''\rho_m(z')(w', w')\geq \frac{L}{\max_{z_1\in F(\overline{W}\cap\overline{V})}\|F'(F^{-1}(z_1))\|^2}\|w'\|^2,\\
 ii')&M'\cdot \max_{z_1\in F(\overline{W}\cap\overline{V})}\|(F^{-1})'(z_1)\|\geq \|D\rho_m(z')\|\geq \frac{M}{\max_{z_1\in F(\overline{W}\cap\overline{V})}\|F'(F^{-1}(z_1))\|}
\end{align*}
for all $z' \in F(\overline{W}\cap\overline{V})$ and $w'\in\mathbb{C}^n$. It is also easy to check that each family $\{\frac{\partial \rho_m}{\partial x_i}\}$ ($i=1,\cdots, 2n$) is equicontinuous on $F(\overline{W}\cap\overline{V})$. Thus $F(U)$ is strongly pseudoconvex with boundary $F(bU)$.
\end{proof}

Let us finish this section stating an invariance property of $0$-uniformly pseudoconvex domains $U$ under biholomorphisms whose domain contain $\overline{U}$; its proof is an easy consequence of the corresponding result for $1$-uniformly pseudoconvex domains. As with $1$-uniformly pseudoconvex domains, we have extended our family of examples of $0$-uniformly pseudoconvex domains. Moreover, we obtain that strongly pseudoconvex domains are invariant under affine isomorphisms between Banach spaces.

\begin{proposition}
The biholomorphic image of a $0$-uniformly pseudoconvex domain $U$ in $\mathbb{C}^n$ is still $0$-uniformly pseudoconvex, when the domain of the biholomorphism contains $\overline{U}$.
\end{proposition}

\section{Cauchy-Riemann equations}\label{S4}

In this final section we will see some solutions to the (inhomogeneous) Cauchy-Riemann equations in infinite dimension. The following notions can be found in \cite[Ch.~V]{M}.

\medskip

\begin{definition}If $X$ and $Y$ denote complex Banach spaces, let $X_{\mathbb{R}}$ and $Y_{\mathbb{R}}$ denote the respective previous spaces seen as real Banach spaces. For every $m \in \mathbb{N}$, $L(^m X_{\mathbb{R}}, Y_{\mathbb{R}})$ denotes the continuous $m$-linear mappings $A:X_{\mathbb{R}}^m \to Y_{\mathbb{R}}$, while $L^a(^m X_{\mathbb{R}}, Y_{\mathbb{R}})$ denotes the continuous $m$-linear mappings $A:X_{\mathbb{R}}^m \to Y_{\mathbb{R}}$ that are \textit{alternating}, i.e.
$$A(x_{\sigma(1)}, \cdots, x_{\sigma(m)})=(-1)^{\sigma} A(x_1, \cdots, x_m), \; \text{ for all } \sigma \in S_m \text{ and } x_1, \cdots, x_m \in X.$$
Also, given $m \in \mathbb{N}$ and $p,q \in \mathbb{N}_0$ such that $p+q=m$, $L^a(^{p,q} X_{\mathbb{R}}, Y_{\mathbb{R}})$ is the subspace of $A \in L^a(^m X_{\mathbb{R}}, Y_{\mathbb{R}})$ such that
$$A(\lambda x_1, \cdots, \lambda x_m)=\lambda^p \bar{\lambda}^q A(x_1, \cdots, x_m), \; \text{ for all } \lambda \in \mathbb{C} \text{ and } x_1, \cdots, x_m \in X;$$
while $L^{apq}(^m X_{\mathbb{R}}, Y_{\mathbb{R}})$ denotes the subspace of all $A \in L(^m X_{\mathbb{R}}, Y_{\mathbb{R}})$ which are alternating in the first p variables and are alternating in the last q variables.

\medskip

Given $A \in L(^m X_{\mathbb{R}}, Y_{\mathbb{R}})$, define its alternating part $A^a \in L^a(^m X_{\mathbb{R}}, Y_{\mathbb{R}})$ by
$$A^a(x_1, \cdots, x_m)=\frac{1}{m!}\sum_{\sigma \in S_m} (-1)^{\sigma} A(x_{\sigma(1)}, \cdots, x_{\sigma(m)}), \; \text{ for all }x_1,\cdots, x_m \in X.$$

\smallskip

Given $U$ an open subset of $X$ and $p,q \in \mathbb{N}_0$, let $C_{p,q}^{\infty}(U,Y):= C^{\infty}(U, L^a(^{p,q} X_{\mathbb{R}}, Y_{\mathbb{R}}))$. Then for each $f \in C_{p,q}^{\infty}(U,Y)$, $\bar{\partial}f \in C_{p,q+1}^{\infty}(U,Y)$ is given by
$$\bar{\partial}f(x)=(m+1)[D''f(x)]^a, \; \text{ for all }x \in U.$$
\end{definition}
\begin{remark} \rm{Since $D''f(x) \in L(X_{\mathbb{R}}, L^a(^m X_{\mathbb{R}}, Y_{\mathbb{R}}) )=L^{a1m}(^{m+1} X_{\mathbb{R}}, Y_{\mathbb{R}})$, then Proposition 18.6 in \cite{M} implies that $\text{ for all }t_1, \cdots, t_{m+1} \in X$, 
\begin{align*}
\bar{\partial}f(x)(t_1, \cdots, t_{m+1})&=(m+1)[D''f(x)]^{a1m}(t_1, \cdots, t_{m+1})\\
&=\frac{m+1}{m+1} \sum_{\sigma \in S_{1m}} (-1)^{\sigma} D''f(x) (t_{\sigma(1)}, \cdots, t_{\sigma(m+1)})
\end{align*}
where $S_{1m}$ denotes the set of all permutations $\sigma \in S_{m+1}$ such that $\sigma(2)< \cdots< \sigma(m+1)$, so
$$\bar{\partial}f(x)(t_1, \cdots, t_{m+1})=\sum_{j=1}^{m+1}(-1)^{j-1} D''f(x)(t_j) (t_1, \cdots, t_{j-1}, t_{j+1}, \cdots, t_{m+1}).$$
}
\end{remark}
The $\bar{\partial}$ problem for $g \in C_{p,q+1}^{\infty}(U,Y)$ asks whether the equation $\bar{\partial}g=0$ implies the existence of $f \in C_{p,q}^{\infty}(U,Y)$ such that $\bar{\partial}f=g$.

\medskip

We are interested in the case when $p=q=0$, since in this case we have the following known result \cite{L}:

\begin{theorem}\label{thm1}
Suppose that $g$ is a complex-valued $(0,1)$-form on $B_{\ell_1}$ that is Lipschitz continuous on all balls $r B_{\ell_1}$, $0<r<1$. If $g$ is a $\overline{\partial}$-closed form then there is a continuously differentiable function $f$ on $B_{\ell_1}$ that solves $\overline{\partial}f=g$. If, in addition, $g$ is $m$ times continuously differentiable, $m=2,3,\dots$, then so is $f$.
\end{theorem}

Let us draw some conclusions from such result. For the following proposition, keep in mind that examples of spaces isomorphic to $\ell_1$ include the infinite-dimensional complemented subspaces of $\ell_1$ \cite{P}.

\begin{proposition}\label{pr1}
Suppose that $U=T(B_{\ell_1})$ for $T=a+L$ an affine isomorphism of Banach spaces $T: \ell_1\to X$. If $g\in C^{\infty}_{0,1}(U, \mathbb{C})$, $\overline{\partial}g=0$, and $g$ is Lipschitz continuous on all $a+r(U-a)$, $0<r<1$, then there exists $f\in C^{\infty}(U,\mathbb{C})$ such that $\overline{\partial}f=g$. 
\end{proposition}
\begin{proof}
Define $\tilde{g}\in C^{\infty}_{0,1}(B_{\ell_1}, \mathbb{C})$ by $\tilde{g}(x)(t)=g(Tx)(Lt)$. Note that $\tilde{g}$ is Lipschitz continuous on all balls $rB_{\ell_1}$, $0<r<1$, because given $r\in (0,1)$ as well as $x_1, x_2\in r B_{\ell_1}$ and $t\in \ell_1$ we have that $Tx_1, Tx_2\in a+r(U-a)$ hence
\begin{align*}
|\tilde{g}(x_1)(t)-\tilde{g}(x_2)(t)|&= |g(Tx_1)(Lt)-g(Tx_2)(Lt)|\\
&\leq M_r \|Tx_1-Tx_2\| \|Lt\|\\
&= M_r\|Lx_1-Lx_2\| \|Lt\|\\
&\leq M_r \|L\|^2 \|x_1-x_2\| \|t\|
\end{align*}
i.e. $\|\tilde{g}(x_1)-\tilde{g}(x_2)\|\leq M_r\|L\|^2 \|x_1-x_2\|$.

\smallskip

Moreover, we also know that $\overline{\partial}\tilde{g}=0$ because given $x\in B_{\ell_1}$ as well as $t_1, t_2 \in \ell_1$ we have that, since $L$ is linear,
\begin{align*}
\overline{\partial}(\tilde{g})(x)(t_1,t_2)&=D''(\tilde{g})(x)(t_1)(t_2)-D''(\tilde{g})(x)(t_2)(t_1)\\
&=D''g(Tx)(Lt_1)(Lt_2)-D''g(Tx)(Lt_2)(Lt_1)\\
&=\overline{\partial}g(Tx)(Lt_1,Lt_2)\\
&=0.
\end{align*}
Then there exists $\tilde{f}\in C^{\infty}(B_{\ell_1}, \mathbb{C})$ such that $\overline{\partial}\tilde{f}=\tilde{g}$. Consequently, $f=\tilde{f}\circ T^{-1}\in C^{\infty}(U, \mathbb{C})$ satisfies that $\overline{\partial}f=g$ since for all $x\in U$ and $t\in X$,
\begin{align*}
\overline{\partial}f(x)(t)&=D''(\tilde{f}\circ T^{-1})(x)(t)\\
&=D''\tilde{f}(T^{-1}x)(L^{-1}t)\\
&=\overline{\partial}\tilde{f}(T^{-1}x)(L^{-1}t)\\
&=\tilde{g}(T^{-1}x)(L^{-1}t)\\
&=g(TT^{-1}x)(LL^{-1}t)\\
&=g(x)(t).
\end{align*}
\end{proof}

Actually, we can extend Proposition \ref{pr1} to images over bounded biholomorphisms at the price of global Lipschitz continuity. Note, however, that we do not know which domains are biholomorphic to the unit ball of $\ell_1$, and that when we deal with automorphisms of the ball of $\ell_1$, these are simply affine isomorphisms \cite{BKU}.

\begin{proposition}\label{pr2}
Suppose that $U=\phi(B_{\ell_1})$ for $\phi:B_{\ell_1}\to U\subset X$ a bounded biholomorphism. If $g\in C^{\infty}_{0,1}(U, \mathbb{C})$, $\overline{\partial}g=0$, and $g$ is Lipschitz continuous on $U$, then there exists $f\in C^{\infty}(U,\mathbb{C})$ such that $\overline{\partial}f=g$. 
\end{proposition}
\begin{proof}
Define $\tilde{g}\in C^{\infty}_{0,1}(B_{\ell_1}, \mathbb{C})$ by $\tilde{g}(x)(t)=g(\phi x)(D'\phi(x) t)$. Then $\tilde{g}$ is Lipschitz continuous on all balls $rB_{\ell_1}$ because so is $\phi$ due to Schwarz' lemma \cite[\S 7]{M} and the convexity of each $rB_{\ell_1}$, also because $D'\phi$ is bounded on each $rB_{\ell_1}$ due to Cauchy Inequality \cite[\S 7]{M} and because from Cauchy integral formula \cite[\S 7]{M} we further have that $D'\phi$ is also Lipschitz on each ball $rB_{\ell_1}$; indeed, we have that for all $r\in (0,1)$ as well as $x_1, x_2 \in r B_{\ell_1}$ and $t\in \ell_1$,
\begin{align*}
|\tilde{g}(x_1)(t)-\tilde{g}(x_2)(t)|&= |g(\phi x_1)(D'\phi(x_1)t)-g(\phi x_2)(D'\phi(x_2)t)|\\
&\leq |g(\phi x_1)(D'\phi(x_1)t)-g(\phi x_2)(D'\phi(x_1)t)|\\&+|g(\phi x_2)(D'\phi(x_1)t)-g(\phi x_2)(D'\phi(x_2)t)|\\
&\leq M \|\phi x_1-\phi x_2\| \|D'\phi (x_1) t\|\\
&+\|g(\phi x_2)\| \|D'\phi(x_1)(t)-D'\phi(x_2)(t)\|\\
&\leq M\cdot \frac{2\|\phi\|}{1-r} \|x_1-x_2\| \cdot \|D'\phi(x_1)\| \|t\|\\
&+(M\cdot 2\|\phi\|+\|g(\phi (0))\|)\frac{\|t\|}{(1-r)/2}\cdot 2\|\phi\| \|x_1-x_2\| \\
&\leq M\cdot \frac{2\|\phi\|}{1-r} \cdot \frac{\|\phi\|}{1-r}\|x_1-x_2\| \|t\|\\
&+(M\cdot 2\|\phi\|+\|g(\phi (0))\|)\frac{4\|\phi\|}{1-r} \|x_1-x_2\| \|t\|.
\end{align*}
Moreover, we also know that $\overline{\partial}\tilde{g}=0$ because given $x\in B_{\ell_1}$ as well as $t_1, t_2 \in \ell_1$ we have that, since $D'\phi(x)$ is linear,
\begin{align*}
\overline{\partial}(\tilde{g})(x)(t_1,t_2)&=D''(\tilde{g})(x)(t_1)(t_2)-D''(\tilde{g})(x)(t_2)(t_1)\\
&=D''g(\phi x)(D'\phi(x)t_1)(D'\phi(x)t_2)-D''g(\phi x)(D'\phi(x)t_2)(D'\phi(x)t_1)\\
&=\overline{\partial}g(\phi x)(D'\phi(x)t_1,D'\phi(x)t_2)\\
&=0.
\end{align*}
Consequently there exists $\tilde{f}\in C^{\infty}(B_{\ell_1}, \mathbb{C})$ such that $\overline{\partial}\tilde{f}=\tilde{g}$. Then $f=\tilde{f}\circ \phi^{-1}\in C^{\infty}_{0,1}(U, \mathbb{C})$ satisfies that, for $x\in U$ and $t\in X$,
\begin{align*}
\overline{\partial}f(x)(t)&=D''(\tilde{f}\circ \phi^{-1})(x)(t)\\
&=D''\tilde{f}(\phi^{-1}x)(D'\phi^{-1}(x)t)\\
&=\overline{\partial}\tilde{f}(\phi^{-1}x)(D'\phi^{-1}(x)t)\\
&=\tilde{g}(\phi^{-1}x)(D'\phi^{-1}(x)t)\\
&=g(\phi\phi^{-1}x)( D'\phi(\phi^{-1}x)D'\phi^{-1}(x)t)\\
&=g(x)(t).
\end{align*}
i.e. $\overline{\partial}f=g$, as we wanted.
\end{proof}

Let us also prove a generalization of the following result in \cite{L0}. 

\begin{theorem}\label{thm2}
If $\Omega\subset \ell_1$ is pseudoconvex and $g\in C_{0,1}(\Omega, \mathbb{C})$ is Lipschitz continuous and $\overline{\partial}$-closed then the equation $\overline{\partial}f=g$ has a solution $f\in C^1(\Omega, \mathbb{C})$.
\end{theorem}

As already mentioned in \cite{L0}, Lipschitz continuous $(0,1)$-forms are not holomorphically invariant. The best we achieved for biholomorphic changes is the following.

\begin{proposition}
Suppose that $U=\phi(\Omega)$ for $\phi:V\subset \ell_1 \to W\subset X$ a bounded biholomorphism such that $\Omega\subset V_{\delta}:=\{z\in V: d(z, bV)>\delta\}$ and that $\Omega$ is convex and open. If $g\in C^{\infty}_{0,1}(U, \mathbb{C})$, $\overline{\partial}g=0$, and $g$ is Lipschitz continuous on $U$, then there exists $f\in C^{\infty}(U,\mathbb{C})$ such that $\overline{\partial}f=g$. 
\end{proposition}
\begin{proof}
Proceed as in the proof of Proposition \ref{pr2}, using Theorem \ref{thm2} instead of Theorem \ref{thm1}, and restricting $\phi$ as well as $D'\phi$ to $\Omega$ instead of $rB_{\ell_1}$.
\end{proof}

\subsection*{Acknowledgements}
The author thanks William B. Johnson, Harold P. Boas, Maite Fern\'andez Unzueta, C\'esar Octavio P\'erez Regalado and Veronique Fischer for several helpful discussions. I also wish to thank Xavier G\'omez-Mont \'Avalos, Eduardo Santill\'an Zer\'on and Loredana Lanzani for helpful suggestions to improve the material in this article. And I thank Isidro Humberto Munive Lima for lendind an ear to the ideas behind Section \ref{S4} of this article.


\end{document}